\documentclass[12pt]{article}
\usepackage{amsfonts, amsmath, amsthm, amssymb,stmaryrd}

\numberwithin{equation}{section}
\newtheorem{theorem}[equation]{Theorem}
\newtheorem{lemma}[equation]{Lemma}
\newtheorem{cor}[equation]{Corollary}
\newtheorem{prop}[equation]{Proposition}
\theoremstyle{definition}

\newtheorem{remark}[equation]{Remark}

\def\QQ{\mathbb{Q}}
\def\ZZ{\mathbb{Z}}
\def\calO{\mathcal{O}}
\def\calR{\mathcal{R}}

\def\fp{\frac{1}{p}}

\def\be{\mathbf{e}}
\def\bv{\mathbf{v}}
\def\bw{\mathbf{w}}

\def\Gcon{\Gamma_{c}}

\DeclareMathOperator{\alge}{alg}

\DeclareMathOperator{\con}{con}
\DeclareMathOperator{\Frac}{Frac}
\DeclareMathOperator{\Gal}{Gal}
\DeclareMathOperator{\Hom}{Hom}
\DeclareMathOperator{\ord}{ord}
\DeclareMathOperator{\perf}{perf}
\DeclareMathOperator{\rank}{rank}
\DeclareMathOperator{\sep}{sep}

\newcounter{fixmectr}

\begin{document}

\title{Frobenius modules and de Jong's theorem}
\author{Kiran S. Kedlaya \\ Department of Mathematics \\ Massachusetts
Institute of Technology \\ 77 Massachusetts Avenue \\
Cambridge, MA 02139 \\
\texttt{kedlaya@math.mit.edu}}
\date{version of December 12, 2004}

\maketitle

\begin{abstract}
Let $k$ be a field of characteristic $p>0$.
A theorem of de~Jong shows that morphisms of modules
over $W(k) \llbracket t \rrbracket$ with Frobenius and connection
structure descend from the completion of $W(k)((t))$. 
A careful reading of de~Jong's proof suggests the possibility
that an analogous theorem holds for modules with only a Frobenius structure.
We show that this analogue holds in one natural
formulation, but fails in a stronger formulation in which
$W(k) \llbracket t \rrbracket$ is replaced by
$W(k \llbracket t \rrbracket)$.
\end{abstract}

\section{Introduction: de Jong's theorem}

In the course of proving the equal characteristic analogue of
Tate's extension theorem for $p$-divisible groups over local fields,
de~Jong \cite{bib:dej4} proved a theorem to the effect that
morphisms of certain $F$-crystals descend. Such crystals can be
regarded as modules equipped with Frobenius and connection structures,
and de Jong's proof primarily uses the Frobenius structure, as 
observed in \cite{bib:methesis} and \cite{bib:mz}. It thus makes 
sense to ask whether one can give a version of de Jong's theorem
using only the Frobenius structure. 
The purpose of this note is to
describe one successful generalization of de Jong's theorem
(Theorem~\ref{thm:main2} below), and to present a counterexample against
a second, more optimistic generalization
(Proposition~\ref{prop:counter}).
In the process, we give a partial exposition of de Jong's theorem which we hope
may be of some value in its own right.

The successful generalization may have some applications
in the theory of $p$-adic differential equations, and consequently
in rigid cohomology. However, it is more likely to be relevant in contexts
where Frobenius structures arise on their own, or in conjunction
with connections that have ``too many singularities''.
Two such contexts are: the work of Berger \cite{bib:berger}
relating $p$-adic Galois representations
to modules with Frobenius and connection (in which
any representation corresponds to a Frobenius structure, but the connection
has poles which can only be removed if the representation is 
de Rham-admissible); and 
the work of Andr\'e and di Vizio \cite{bib:andre divizio} on
a $q$-difference analogue of the theory of $p$-adic differential equations
(in which the form of the connection is perturbed but the notion of a Frobenius
structure does not change).
 
Before proceeding further,
we recall the statement of de~Jong's theorem
\cite[Theorem~9.1]{bib:dej4}.
Let $k$ be an algebraically closed field of characteristic $p>0$,
and let $W$ denote the ring of ($p$-typical) Witt vectors over $k$.
Put $\Omega = W \llbracket t \rrbracket$,
and let $\Gamma$ be the $p$-adic completion of $W((t))$;
elements of $\Gamma$ can be written as power series $\sum_{n \in \ZZ}
c_n t^n$ with $c_n \in W$ and $c_n \to 0$ as $n \to -\infty$.
Fix a power $q$ of $p$,
let $\sigma: W \to W$ denote the $(\log_p q)$-th power of the
Witt vector Frobenius, and
extend $\sigma$ to a map on $\Gamma$
sending $\sum c_n t^n$ to $\sum c_n^\sigma t^{nq}$.

An \emph{$F$-crystal} over $\Omega$ is a finite free
$\Omega$-module $M$ equipped with a $W$-linear connection
$\nabla: M \to M \otimes dt$ (i.e., a map satisfying the Leibniz rule
$\nabla(r\bv) = r\nabla(\bv) + \frac{dr}{dt} \bv \otimes dt$ for
$r \in \Omega$ and $\bv \in M$) and an isogeny $F: \sigma^* M \to M$
of $\Omega$-modules with connection, called its ``Frobenius structure''.
(An \emph{isogeny} here means a morphism
whose kernel and cokernel are killed by some power of $p$.
One could turn isogenies into isomorphisms by working over
$\Omega[\frac 1p]$ instead of $\Omega$, as in
\cite{bib:mz}; we will not do so here.)
Then $F$ induces a $\sigma$-linear
map from $M$ to itself, which we denote by the same letter $F$.

In this language, de Jong's theorem is the following. (Strictly
speaking, de Jong only addressed the case $q=p$, but the general case
follows immediately by a ``restriction of scalars'' argument, as in
\cite[Proposition~6.11]{bib:mecrew}.)
\begin{theorem}[de Jong] \label{thm:main1}
Let $M$ be an $F$-crystal over $\Omega$, and suppose
$\bv \in M \otimes \Gamma$ satisfies $\nabla \bv = 0$
and $F \bv = p^\ell \bv$ for some
integer $\ell$. Then $\bv \in M$.
\end{theorem}
Note that the category of $F$-crystals contains internal Homs
(up to Tate twists),
so this implies that morphisms over $\Gamma$ between $F$-crystals over
$\Omega$ are themselves defined over $\Omega$. (It is in this form in which
Theorem~\ref{thm:main1} implies the extension theorem for $p$-divisible
groups, by application to the corresponding Dieudonn\'e modules.)

It should be noted that the choice of $\sigma$ above was a bit artificial,
made for convenience in the proof. Define a \emph{Frobenius lift}
on $\Omega$ to be any 
ring map $\sigma: \Omega
\to \Omega$ extending the $(\log_p q)$-th power of the
Witt vector Frobenius on $W$ and lifting the
$q$-power map modulo $p$; we say $\sigma$ is \emph{standard} if
$t^\sigma = t^q$.
Then as described in \cite[2.4]{bib:katz},
the data of an $F$-crystal $M$ over $\Omega$ determines a Frobenius structure
$F_\sigma: \sigma^* M \to M$ for any Frobenius lift $\sigma$ on $\Omega$,
so de~Jong's theorem is equally true for any $\sigma$.

In this terminology, the purpose of this note is to explore the
possibility of establishing Theorem~\ref{thm:main1} in the presence
of a Frobenius structure but not a connection.
Our main affirmative result in this direction is the following.
Define an \emph{$F$-module} over $\Omega$, with respect to the Frobenius
lift $\sigma$, to be a finitely generated torsion-free
$\Omega$-module $M$ equipped with an isogeny
$F: \sigma^* M \to M$
of $\Omega$-modules. (Note that in the absence of a connection,
the category of $F$-modules depends on the choice of $\sigma$.
Note also that the freeness hypothesis in the theorem
is needed for trivial reasons; see Remark~\ref{R:need free}.)
\begin{theorem} \label{thm:main2}
Let $M$ be a free $F$-module over $\Omega$
with respect to a Frobenius lift $\sigma$,
and suppose $\bv \in M \otimes \Gamma$
satisfies $F\bv = p^\ell \bv$ for some integer $\ell$. Then $\bv \in M$. 
\end{theorem}
Most of the work in proving this theorem goes into an auxiliary
result extending \cite[Proposition~6.4]{bib:dej4}, which may be
of interest in its own right. That result is Theorem~\ref{thm:64},
which asserts that
given an $F$-module $M$ over $\Omega$, any submodule of $M \otimes
\Gcon$, where $\Gcon$ is the ``overconvergent'' subring of $\Gamma$
(see Section~\ref{sec:sketch}), is actually defined over $\Omega$.

We conclude this introduction by describing our negative observation.
An $F$-module over $\Omega$ gives rise naturally to an $F$-module over
the Witt ring $W(k\llbracket t \rrbracket)$, so one can ask
even more generally whether the analogue of Theorem~\ref{thm:main2}
holds for $F$-modules over $W(k \llbracket t \rrbracket)$ (namely,
whether morphisms descend to this ring from $W(k((t)))$). This
analogue turns out to fail; we give a counterexample in
Proposition~\ref{prop:counter}.

\subsection*{Acknowledgments}
Thanks to Thomas Zink for helpful discussions, notably in suggesting the
Witt vector generalization and in providing a copy of \cite{bib:mz}.
Thanks also to the Institute for Advanced Study for its hospitality.
The author was supported by the National Science Foundation (grant
number DMS-0111298).

\section{Initial observations}

In this section, we make some initial observations in connection with
Theorem~\ref{thm:main2}. These are based on the fact that
(as is evident upon reading \cite{bib:dej4} closely)
the arguments of \cite{bib:dej4} use the connection data of
an $F$-crystal only in an incidental fashion. 
For the ``standard'' Frobenius $t \mapsto t^q$, this was observed in
\cite[Theorem~6.3.1]{bib:methesis} and \cite[Theorem~3]{bib:mz}
(though both arguments include similar subtle errors; see 
Remark~\ref{rem:error}).
By being more careful, one can also eliminate the dependence on the
standard Frobenius (and the erroneous arguments).

To begin with, note that it is not necessary to work with completely
general Frobenius lifts.
We call a Frobenius lift $\sigma$ \emph{zero-centered} if $t^\sigma$
is divisible by $t$. Then one has the following easy result.
\begin{lemma} \label{lem:centered}
Every Frobenius lift is conjugate (by an automorphism of $\Omega$
reducing to the identity map modulo $p$) to a zero-centered
Frobenius lift.
\end{lemma}
\begin{proof}
Given the Frobenius lift $\sigma$, put $t^\sigma = \sum_{i=0}^\infty a_i t^i$.
Then the map $\tau: pW \to pW$ defined by
\[
\tau(x) = \left( \sum_{i=0}^\infty a_i x^i \right)^{\sigma^{-1}}
\]
is a contraction mapping, since
\[
|\tau(x) - \tau(y)| \leq \max_{i>0} \{ |a_i| \cdot |y-x| \cdot 
|(y^i-x^i)/(y-x)| \}
< |y-x|.
\]
(The latter holds because $|a_i| < 1$ if $i \neq q$, and $|y^q - x^q| <
|y-x|$.)
Thus $\tau$ has a unique fixed point $c$, and the equation $\tau(c) = c$
means that $t^\sigma - c^\sigma$ is divisible by $t-c$. Thus the
conjugation by $t \mapsto t-c$ does the job.
\end{proof}

In practice one usually wants to work only with $F$-modules which are
free. A simple observation, adapted from \cite[Lemma~6.1]{bib:dej4},
makes this possible.
\begin{lemma} \label{lem:free}
Let $M$ be an $F$-module over $\Omega$. Then $M' = 
\Hom_\Omega(\Hom_\Omega(M,\Omega))$ is free and admits a natural
$F$-module structure under which the natural map $M \to M'$ is an isogeny
of $F$-modules.
\end{lemma}
\begin{proof}
The freeness of $M'$ follows from the fact that $\Omega$ is a regular local
ring (of dimension 2), so the dual of any finitely generated $\Omega$-module
is free. (Thus $\Hom_\Omega(M,\Omega)$ is free, and thence $M'$ is free.)
The natural map $M \to M'$ is injective and has finite length cokernel
(so in particular is an isogeny). Thus we may choose a positive integer $n$
such that $t^n \bv' \in M$ and $p^n \bv' \in M$ whenever $\bv' \in M'$.

For $\bv' \in M'$, we have $(t^\sigma)^n F(p^n \bv') = p^n F(t^n \bv')$;
since $M'$ is free and $\Omega$ is factorial, 
there must exist $\bw \in M'$ such that
$F(t^n \bv') = (t^\sigma)^n \bw$ and $F(p^n \bv') = p^n \bw$. We set
$F(\bv') = \bw$; this map clearly has the desired properties.
\end{proof}

We will need the classification of $F$-modules of rank 1, following
\cite[Lemma~6.2]{bib:dej4}.
\begin{lemma} \label{lem:rank1}
Let $M$ be an $F$-module of rank $1$ over $\Omega$. Then $M$ is isogenous
to $\Omega$ equipped with $p^\ell \sigma$ for some nonnegative integer $\ell$;
if $M$ is free, then the isogeny can be taken to be an isomorphism.
\end{lemma}
\begin{proof}
By Lemma~\ref{lem:free}, we may reduce to the case where $M$ is free.
Choose a generator $\bv$ of $M$ and write $F\bv = p^{\ell} c\bv$ for $c \in \Omega$
not divisible by $p$. Recall that the definition of an $F$-module requires
that $F$ be an isogeny, so that $\Omega/c\Omega$ is annihilated by a power of $p$.
That is only possible if $c$ is a unit in $\Omega$, that is,
$c \not\equiv 0 \pmod{(p,t)}$.

Since $k$ is algebraically closed, we can find $a \in W^*$ such that
$b = a^{-1} c a^\sigma \equiv 1 \pmod{(p,t)}$.
Now note that the infinite product $b b^\sigma b^{\sigma^2} \cdots$ converges
$(p,t)$-adically to a limit $u \in \Omega^*$ satisfying
$bu^\sigma = u$. If we put $\bw = au \bv$, we then have
$F\bw = p^\ell \bw$, as desired.
\end{proof}
\begin{remark} \label{R:need free}
Note that the freeness hypothesis in the last assertion
of Lemma~\ref{lem:rank1} is needed, else we could take
$M$ to be the ideal $(p,t)$ in $\Omega$ equipped with the standard
Frobenius lift, which is isogenous to $\Omega$ via the natural inclusion
but is not isomorphic to $\Omega$. Indeed, this same example shows that
the freeness hypothesis is needed in Theorem~\ref{thm:main2}, since
$\bv = 1 \in M \otimes \Gamma \cong \Gamma$ does not lie in $M$.
\end{remark}

To conclude this section, we state a dual form of
Theorem~\ref{thm:main2}, which is what we
will actually prove.
\begin{theorem} \label{thm:main3}
Let $M$ be an $F$-module over $\Omega$ with respect to a Frobenius lift
$\sigma$. Suppose $\phi: M \to \Gamma$ is an $\Omega$-linear map with
the property that for some integer $\ell \geq 0$,
\[
\phi(F\bv) = p^{\ell} \phi(\bv)^\sigma \qquad \mbox{for all $\bv \in M$}.
\]
Then $\phi(M) \subseteq \Omega$.
\end{theorem}
To reduce Theorem~\ref{thm:main2} to Theorem~\ref{thm:main3},
take $M, \bv$ as in Theorem~\ref{thm:main2} and put $d = \rank M$. 
Then exterior product with $\bv$ induces a map of $F$-modules
from $\wedge^{d-1} M$, viewed as an $F$-module in the natural fashion, to 
$(\wedge^d M) \otimes \Gamma$
with its natural $F$-action multiplied by $p^\ell$.
Since $\wedge^d M$ has rank 1 over $\Omega$, Lemma~\ref{lem:rank1} implies
that $\wedge^d M$ is isogenous to $\Omega$ with Frobenius given by $p^m \sigma$
for some integer $m$. Thus $\bv$ induces a map $\phi: \wedge^{d-1} M \to
\Gamma$ such that $\phi(F\bv) = p^m \phi(\bv)^\sigma$.
Theorem~\ref{thm:main3} implies that $\phi(\wedge^{d-1} M) \subseteq \Omega$,
so $\bv$ belongs to the double dual of $M$ over $\Omega$. Since $M$
is free, that double dual coincides with $M$, so $\bv \in M$
as desired.

\section{Sketch of the proof}
\label{sec:sketch}

The proof of Theorem~\ref{thm:main3} runs parallel to that of
\cite[Theorem~9.1]{bib:dej4}. In this section, we collect some statements
analogous to various statements in \cite{bib:dej4}, then describe how
they are used to prove Theorem~\ref{thm:main3}. It will then remain to
establish some of the analogues in our more general setting.

Before proceeding further, we recall some 
auxiliary rings from \cite{bib:dej4}.
To begin with, let $\Gcon$
denote the subring of $\Gamma$ consisting of
series $\sum c_n t^n$ such that $v(c_n) + rn \to \infty$ as $n \to -\infty$
for some $r>0$ depending on the series. (These are sometimes
called ``overconvergent''
elements of $\Gamma$.)
The discrete valuation ring $\Gcon$ is not complete; it is henselian
(see, e.g., \cite[Lemma~3.9]{bib:mecrew}), but we won't need this fact.
Put $\Gamma_2 = W(k((t))^{\alge})$; we may
embed $\Gamma$ into $\Gamma_2$
by first identifying the completion
of the direct limit $\Gamma \stackrel{\sigma}{\to} 
\Gamma \stackrel{\sigma}{\to} \cdots$
with $W(k((t))^{\perf})$, then using Witt vector functoriality to embed
$W(k((t))^{\perf})$ into $W(k((t))^{\alge})$.
We then have a subring $\Gamma_{2,c}$ of $\Gamma_2$ analogous to
$\Gamma_c$, consisting of series $\sum p^i [c_i]$
such that $i + r v_t(c_i) \to \infty$ as $i \to \infty$ for some
$r>0$. (Here $v_t$ is the $t$-adic valuation on $k((t))^{\alge}$
and brackets denote Teichm\"uller lifts.) Note that
$\Gamma \cap \Gamma_{2,c} = \Gamma_c$.
(In the notation of \cite{bib:mecrew}, the rings
$\Gcon, \Gamma_{2}, \Gamma_{2,c}$ would be denoted $\Gamma_{\con},
\Gamma^{\alge}, \Gamma^{\alge}_{\con}$, respectively.)
For $R$ one of the aforementioned rings, or indeed
any ring equipped with an endomorphism $\sigma$,
we define an \emph{$F$-module} over $R$, with respect to $\sigma$,
as a finite torsion-free $R$-module $M$ equipped with an isogeny
$F: \sigma^* M \to M$.

To describe a key result of \cite[Section~5]{bib:dej4} we will use, we must
recall the notion of slopes and the Dieudonn\'e-Manin classification of
$F$-modules over $R$, for $R$ a complete discrete valuation ring of mixed
characteristics $(0,p)$ with algebraically closed residue field
\cite{bib:manin}, \cite{bib:katz}.
Namely, for any $F$-module $M$ over $R$, there exists a positive integer $m$
and a basis $\bv_1, \dots, \bv_n$ of $R$ such that
$F^m \bv_i = p^{\ell_i} \bv_i$ for some nonnegative integers $\ell_i$.
The multiset $\{\ell_1/m, \dots, \ell_n/m\}$ does not depend on any choices;
its elements are called the \emph{slopes} of $M$.
For $M$ an $F$-module over $\Gamma$ (or $\Omega$ or $\Gcon$), 
we define slopes by base extension to $\Gamma_2$.

One consequence of \cite[Corollary~5.7]{bib:dej4}
 is the following.
\begin{prop} \label{prop:57}
Let $M$ be an $F$-module over $\Gamma_{2,c}[p^{1/b}]$ for some
positive integer $b$, and suppose
$\phi: M \to \Gamma_{2}[p^{1/b}]$ is a $\Gamma_{2,c}[p^{1/b}]$-linear 
map such that
for some nonnegative integer $\ell$, we have $\phi(F\bv)
= p^{\ell/b} \phi(\bv)^\sigma$. If $\ell/b$ is greater than every slope
of $M$, then $\phi = 0$.
\end{prop}
\begin{proof}
In case $q=p$ and $\sigma$ is standard, this is part (iii) of
\cite[Corollary~5.7]{bib:dej4}. However,
the Frobenius on $\Gamma_{2,c}$
or $\Gamma_2$ is independent of the initial choice of $\sigma$ (it is the
Witt vector Frobenius), so changing $\sigma$ does not affect the truth
of this result.

The case of $q$ general can be treated by imitating the arguments of
\cite{bib:dej4} \emph{mutatis mutandis}. More explicitly,
one may replace the invocation of \cite[Proposition~5.5]{bib:dej4}
in the proof of \cite[Corollary~5.7]{bib:dej4} by 
\cite[Proposition~5.11]{bib:mecrew}.
This argument remains encumbered by the running hypothesis in \cite{bib:mecrew}
that the Frobenius lift $\sigma$ is a power of a $p$-power Frobenius lift;
however, that assumption is only actually used starting with
\cite[Proposition~6.11]{bib:mecrew}, so no prior results of
\cite{bib:mecrew} actually depend on it.
\end{proof}

A further consequence is the following
result.
\begin{prop} \label{prop:topslope}
Let $M$ be a nonzero $F$-module over $\Gcon$, and suppose
$\phi: M \to \Gamma$ is an injective $\Gcon$-linear map such that
for some nonnegative integer $\ell$, $\phi(F\bv) = p^{\ell}
\phi(\bv)^{\sigma}$
for all $\bv \in M$. Then the largest slope of $M$ is equal to $\ell$
with multiplicity $1$, and $\phi^{-1}(\Gcon)$ is a rank $1$ sub-$F$-module
of $M$ of slope $\ell$.
\end{prop}
\begin{proof}
This result is \cite[Corollary~8.2]{bib:dej4} when the Frobenius lift
is standard and $q=p$. It is \cite[Lemma~4.2]{bib:mefull} if the Frobenius
lift is a power of a $p$-power Frobenius lift; however, the argument
goes through unchanged without this hypothesis (as in the previous proof).
\end{proof}

We quote one additional result, namely
\cite[Proposition~7.1]{bib:dej4}, which allows us to split certain
exact sequences under a restriction on slopes. 
The proof of this result for a general Frobenius
is the same as for the standard Frobenius 
(except that \cite[Lemma~6.1]{bib:dej4}
is replaced by our Lemma~\ref{lem:free}).
\begin{prop} \label{prop:split}
Let $M$ be an $F$-module over $\Omega$. Suppose that $N \subseteq M$
is a saturated rank $1$ sub-$F$-module of slope $\ell$ and that the slopes of 
$M/N$ are all less than $\ell$. Then there is a sub-$F$-module $N'$ of $M$
such that the natural map $N \oplus N' \to M$ is an isogeny.
\end{prop}

Finally, we state one result which we cannot cite directly from
\cite{bib:dej4}; its proof will occupy much of the rest of the paper.
\begin{theorem} \label{thm:64}
Let $M$ be an $F$-module over $\Omega$,
and let $N_c$ be a saturated sub-$F$-module of $M_c = M \otimes \Gcon$.
Then $N_c = N \otimes \Gcon$ for some saturated sub-$F$-module 
$N$ of $M$.
\end{theorem}
\begin{remark}
Note that Theorem~\ref{thm:64} does not hold with $\Gcon$ replaced by
$\Gamma$, even in the presence of a connection; one obtains a counterexample
by considering the middle cohomology of the Legendre family of elliptic curves
near a fibre with supersingular reduction, together with its
``unit root'' sub-$F$-module.
\end{remark}

With these ingredients, we can now give the proof of
Theorem~\ref{thm:main3} conditioned on
Theorem~\ref{thm:64}.
\begin{proof}[Proof of Theorem~\ref{thm:main3}]
We may of course assume that $\phi$ is injective, and in particular
that $M$ is torsion-free.
Given $\phi$, let $\phi_c$ be the composite map
\[
M_c \stackrel{\phi \otimes 1}{\to}
\Gamma \otimes_{\Omega} \Gcon \stackrel{\mu}{\to}
\Gamma \otimes_{\Gcon} \Gcon = \Gamma,
\]
where $\mu$ denotes the multiplication map.

We first proceed assuming that $\phi_c$ is injective.
Put $N_c = \phi_c^{-1}(\Gcon)$;
then by Proposition~\ref{prop:topslope},
$N_c$ is a saturated rank 1 sub-$F$-module
of $M_c$ of slope $\ell$, and all other slopes of $M$ are
strictly less than $\ell$. 
Let $N$ be a saturated sub-$F$-module of $M$ with $N \otimes \Gcon =
N_c$ given by Theorem~\ref{thm:64}; 
then $N$ also has rank 1, so Lemma~\ref{lem:rank1}
implies that $N$ is isogenous as an $F$-module
to $\Omega$ with Frobenius action given by $p^m \sigma$ for some integer $m$.
The slope of the latter is $m$, so we must have $\ell = m$. 
Choose $\bv \in N$ such that $F \bv = p^\ell \bv$; then
$g = \phi(\bv)$ satisfies $g^\sigma = g$, so it must belong to
$\ZZ_p$. Hence $\phi(N) \subseteq \Omega[\fp] \cap \Gamma = \Omega$.

By Proposition~\ref{prop:split}, there is an isogeny
$N \oplus N_1 \to M$ for some saturated sub-$F$-module $N_1$ of $M$.
The map $N_1 \to M \stackrel{\phi}{\to} \Gamma$ is the composite of 
injections, so is also injective. On the other hand,
the slopes of $N_1$ are all strictly less than $\ell$, so the 
composite
\[
N_1 \otimes_{\Omega} \Gamma_{2,c}[p^{1/b}] =
(N_1 \otimes_\Omega \Gamma_c) \otimes_{\Gamma_c} \Gamma_{2,c}[p^{1/b}]
\to (\Gamma \otimes_{\Gamma_c} \Gamma_{2,c})[p^{1/b}]
\stackrel{\mu}{\to} \Gamma_2[p^{1/b}]
\]
is zero by Proposition~\ref{prop:57}.

Since the map $\mu: \Gamma \otimes_{\Gamma_c} \Gamma_{2,c} \to \Gamma_2$
is injective by
\cite[Proposition~4.1]{bib:mefull} (again with the appropriate modifications
in case $\sigma$ is not a power of a $p$-power Frobenius lift), 
the map $N_1 \otimes_{\Omega} \Gamma_c
\to \Gamma$ becomes zero after tensoring over $\Gamma_c$ with $\Gamma_{2,c}$,
and so is itself zero. (Here it matters that $\Gamma_{2,c}$
is flat over $\Gamma_c$, which is clear: $\Gamma_c$ is a principal ideal
domain, so a module over $\Gamma_c$ 
is flat if and only if it is torsion-free, which
$\Gamma_{2,c}$ evidently is.) Since $\phi_c$ is injective, it follows that
$N_1 \otimes_\Omega \Gamma_c = 0$. By Lemma~\ref{lem:free},
there is an isogeny $N_1 \to N_1'$ with $N_1'$ free over $\Omega$;
tensoring with $\Gamma_c$, we see that $N_1' \otimes_\Omega \Gamma_c
\to (N_1'/N_1) \otimes_\Omega \Gamma_c$ is a bijection. Since the source
of this map is free over $\Gamma_c$ and the target is killed by a power of $p$,
both must vanish. Thus $N_1' = 0$; since $N_1$ is isogenous to $N_1'$, it is
also zero. In other words,
the inclusion $N \subseteq M$ is an
isogeny. Since $\phi(N) \subseteq \Omega$, we have
$\phi(M) \subseteq \Omega[\frac{1}{p}] \cap \Gamma = \Omega$, as desired.

We now treat the case where $\phi_c$ is not injective.
Put $N'_c = \ker(\phi_c) \subseteq M_c$.
Clearly $N'_c$ is a saturated sub-$F$-module of $M_c$, so
by Theorem~\ref{thm:64} again,
$N'_c = N' \otimes_\Omega \Gcon$ for some saturated sub-$F$-module $N'$ of $M$,
Clearly $N' \subseteq 
\ker(\phi)$; thus we may apply the previous argument to $M/N'$
to deduce that again $\phi(M) \subseteq \Omega$, as desired.
\end{proof}

\section{Some nonarchimedean function theory: Dwork's trick}

As we have just seen,
to complete the proof of Theorem~\ref{thm:main3}, we must supply
Theorem~\ref{thm:64} to replace \cite[Proposition~6.4]{bib:dej4}.
This is harder than one might expect; as noted in \cite{bib:mz}, this is
perhaps the subtlest aspect of the proof of \cite[Theorem~9.1]{bib:dej4}. 

We begin with a version of the ``Dwork trick'' \cite[Lemma~6.3]{bib:dej4}.
Let $\calR^+$ denote the subring of $W[\frac 1p]\llbracket t \rrbracket$
consisting of series $\sum_{n =0}^\infty c_n t^n$ such that for any $r>0$,
$v(c_n) + rn \to \infty$ as $n \to \infty$.
That is, such series converge for $t$ in the open unit disc in $W^{\alge}$,
i.e., are the rigid analytic functions on the open unit disc $D$.
(Correspondingly, this ring is denoted $\Gamma(D, \calO_D)$ in
\cite{bib:dej4}.)
For each $r>0$, the function
\begin{equation} \label{eq:wr}
w_r\left(\sum c_i t^i \right) = \min_i \{v(c_i) + ri\}
\end{equation}
gives a valuation on $\calR^+$ (which induces
the supremum norm on the disc of radius $p^{-r}$).
The ring $\calR^+$ is complete for the Fr\'echet topology
generated by the valuations $w_r$ over all $r>0$, or even any
cofinal set of $r>0$ (e.g., all rational $r>0$).

The proof of Dwork's trick in \cite{bib:dej4}
uses crucially the property of the standard Frobenius
that $t^\sigma$ is divisible by $t^2$, so that repeatedly applying
$\sigma$ to $t$ gives terms which converge to zero $t$-adically.
To make the argument work more generally, one must be a bit more careful.
\begin{lemma} \label{lem:wr}
Suppose that $\sigma$ is zero-centered.
For any $s>0$ and any $r \in (0, \min\{qs,s+1\}]$,
we have
\[
w_s((t^h)^\sigma) \geq rh
\]
for all $h > 0$.
\end{lemma}
\begin{proof}
Write $t^\sigma = t^q + pu$, where $u \equiv 0 \pmod{t}$ because
$\sigma$ is zero-centered; then $w_s(pu) \geq 1 + s$, and so
\[
w_s((t^h)^\sigma) \geq \min_{0 \leq i \leq h} \{ qsi + (s+1)(h-i) \}
\geq \min\{qsh, (s+1)h\}.
\]
The claim thus follows.
\end{proof}

\begin{prop} \label{prop:dwork}
Let $M$ be a free $F$-module over $\Omega$, with respect to a Frobenius
lift $\sigma$. Then $M \otimes_\Omega \calR^+$
admits a basis $\bv_1, \dots, \bv_n$ such that 
$F^m \bv_i = p^{\ell_i} \bv_i$ for some integers $m>0$ and $\ell_i \geq 0$.
\end{prop}
The proposition can also be formulated for free $F$-modules over
$\calR^+$; the proof goes through unchanged.
\begin{proof}
We may assume $\sigma$ is zero-centered.
By the Dieudonn\'e-Manin classification, we can choose a basis
$\be_1, \dots, \be_n$ of $M \otimes \Omega[\fp]$
such that $F^m \be_i \equiv p^{\ell_i} \be_i \pmod{tM}$
for some integers $m>0$ and $\ell_i \geq 0$.
Define the $n \times n$  matrix $\Phi$ by
\[
F^m\be_j = \sum_i \Phi_{ij} \be_i,
\]
let $I$ denote the $n \times n$ identity matrix,
and let $D$ be the $n \times n$ diagonal matrix with $D_{ii} = p^{\ell_i}$.
Define $c \in pW$ by $t^\sigma \equiv ct \pmod{t^2}$, and choose
$h_0$ large enough so that
\[
h_0 v(c) \geq \max_{i,j} \{v(D_{ii}) - v(D_{jj})\}.
\]

We first verify that for any integer $h>0$, there exists an
invertible matrix $U_h$ over $W[\frac{1}{p}]\llbracket t \rrbracket$
with $U_h \equiv I \pmod{t}$ such that $U_h^{-1} \Phi U_h^{\sigma^m}
\equiv D \pmod{t^h}$. For $h=1$, we may take $U_h = I$. Given $U_h$,
put
\[
Y_h = U_h^{-1} \Phi U_h^{\sigma^m}D^{-1} - I,
\]
and define the matrix $X_h$ over $W[\frac 1p]$ by $Y_h + X_ht^h
\equiv 0 \pmod{t^{h+1}}$.
Then there exists a matrix $V_h$ over $W[\frac 1p]$ such that
\begin{equation}
X_h = c^h D V_h^{\sigma^m} D^{-1} - V_h,
\end{equation}
since this equation reduces to a separate equation
$(X_h)_{ij} = c^h D_{ii} (V_h)_{ij}^{\sigma^m}D_{jj}^{-1} - (V_h)_{ij}$ 
for each
matrix entry
$(V_h)_{ij}$, and each of these is solvable when $k$ is algebraically closed.
Moreover, the solution is unique when $v(c^h D_{ii} D_{jj}^{-1}) \neq 0$, which
is automatic for $h > h_0$.
We may thus take $U_{h+1} = U_h(I + V_h t^h)$.

Because $U_h$ modulo $t^h$ determines $U_{h+1}$ modulo $t^{h+1}$ for
$h > h_0$, we can take $t$-adic limits to obtain an invertible matrix $U$
over $W[\frac 1p]\llbracket t \rrbracket$ such that $U^{-1} \Phi U^{\sigma^m}
= D$. We next show that $w_r(U)$ is defined for $r$ sufficiently large. Here
we say ``$w_r(x)$ is defined'' for $x = \sum c_i t^i \in 
W[\frac{1}{p}]\llbracket t \rrbracket$
if the right side of \eqref{eq:wr} is defined, in which case we use
\eqref{eq:wr} as the definition. Also, we define $w_r$ of a matrix as 
the minimum over its entries.

Choose $r_0>h_0 v(c)$ such that $w_r(Y_{h_0})$
is defined and positive for $r \geq r_0$. 
(This is possible because by induction on $h$,
$w_r(Y_h)$ is defined for $r$ sufficiently large, depending on $h$.)
Note that if $U_{h+1} = U_h(I+V_ht^h)$ as above, then
\[
Y_{h+1}
= (I+V_ht^h)^{-1} Y_h (I + DV_ht^hD^{-1})^{\sigma^m}.
\]
For $h \geq h_0$, we have
$w_r(D_{ii} (t^h)^{\sigma^m} D_{jj}^{-1}) > 
w_r(t^h)$; also, if $w_r(Y_h)$ is defined, then $w_r(V_h t^h) \geq w_r(X_h t^h)
\geq w_r(Y_h)$.
Hence if $w_r(Y_h)$ is defined and positive, then
$w_r(Y_{h+1}) \geq w_r(Y_h)$; by induction on $h$, we see that
$w_r(Y_h)$ is defined and positive for all $h \geq h_0$.
Moreover, we have $w_r(Y_{h+1})
\geq w_r(Y_h)$, with strict inequality when $w_r(Y_h-V_h t^h) > w_r(Y_h)$.
The latter
happens infinitely often unless the $Y_h$ eventually become all zero,
in which case so do the $V_h$.
Hence $w_r(U_h^{-1} U_{h+1} - I) \to \infty$ as $h \to \infty$, so that
$w_r(U)$ is defined.

From $r_0$, define the sequence $r_0, r_1, \dots$ by setting
$r_{i+1} = \max\{r_i/q,r_i-1\}$. By Lemma~\ref{lem:wr}, if $w_r(U)$ is
defined for $r \geq r_i$, then $w_r(U^{\sigma^m})$ is defined for
$r \geq r_{i+m}$. Since $U = \Phi U^{\sigma^m} D^{-1}$, it then follows
that $w_r(U)$ is defined for $r \geq r_{i+m}$. Since $w_r(U)$ is defined for
$r \geq r_0$, we deduce by induction on $i$ that $w_r(U)$ is defined for
$r \geq r_{im}$ for all $i$. Since $r_{im} \to 0$ as $i \to \infty$, it
follows that $w_r(U)$ is defined for all $r>0$; in other words,
$U$ has entries in $\calR^+$.

The same argument as above,
applied to the inverse transpose of $U$, shows that $U^{-1}$
also has entries in $\calR^+$. We may thus define a basis
$\bv_1, \dots, \bv_n$ of $M \otimes \calR^+$ by
\[
\bv_j = \sum_i U_{ij} \be_i
\]
and we will have $F^m\bv_i = p^{\ell_i} \bv_i$, as desired.
\end{proof}

\begin{remark} \label{rem:product}
Note that in the above argument, if $D$ is a scalar matrix,
we may take $h_0 = 0$. In particular, in the one-dimensional case,
one may deduce that if $\sigma$ is zero-centered
and $u \in \Omega[\fp]$ is congruent to 1 modulo
$t$, then the infinite product $uu^\sigma u^{\sigma^2} \cdots$
converges in $\calR^+$.
\end{remark}

\section{More nonarchimedean function theory: descent of submodules}

We now proceed towards a proof of Theorem~\ref{thm:64}.
Even for $\sigma$ standard, the proof
in the absence of a connection is subtle, and the arguments
in \cite{bib:methesis} and \cite{bib:mz} are inadequate on this point;
see Remark~\ref{rem:error}.

By a \emph{principal unit} of $\Omega$, we will mean a unit which,
as a series in $t$, has constant coefficient 1.
\begin{lemma} \label{lem:nonvanish}
Suppose $a \in \Gcon[\fp]$ has the property that the series $a$ and $a^{-1}$
converge for $0 < |t| < 1$. Then $a = c t^n u$ for some $c \in W[\fp]$,
some integer $n$ and some principal unit $u$ of $\Omega$.
\end{lemma}
\begin{proof}
We can choose $c,n,u$ such that $b = a/(c t^n u)$ has the form
$1 + \sum_{i=1}^\infty b_i t^{-i}$ with $|b_i| < 1$ for all $i$, as in 
\cite[Proposition~6.5]{bib:mecrew}. It then remains to prove that $b = 1$,
given that $b$ converges to nonzero values for $0 < |t| < 1$.
Suppose on the contrary that $b \neq 1$. Since $b$ converges for
$0 < |t| < 1$, the sequence $\{|b_i|^{1/i}\}_{i=1}^\infty$ 
tends to zero. In particular,
there is some $i$ which maximizes $|b_i|^{1/i}$; choose the largest
such $i$ and put $\rho = |b_i|^{1/i}$, so that $0 < \rho < 1$. By Weierstrass
preparation, $b$ factors as a nonconstant polynomial in $t^{-1}$
whose roots have
absolute value $\rho^{-1}$, times a series in $t^{-1}$ which is invertible
on the disc $|t^{-1}| \leq \rho^{-1}$. Thus $b$ vanishes
for some $t$ with $0 < |t| < 1$, contradiction. Hence $b=1$
and $a = ct^n u$, as desired.
\end{proof}

Let $D$ denote the open unit disc over $W$, or more precisely,
the set of elements of $W^{\alge}$ of norm less than 1.
Fix a choice of an extension of $\sigma$ to an automorphism of $W^{\alge}$
(which exists by standard Galois theory).
Define $a_0, a_1, \dots, \in W$ by the formula $t^\sigma = \sum a_i t^i$, and
let $\tau: D \to D$ denote the map 
\[
\tau(x) = \left( \sum_{i=0}^\infty a_i x^i \right)^{\sigma^{-1}},
\]
so that if we view $f \in \Omega$ as a function of $t$, then
for $y \in D$,
\[
f^\sigma(y) = f(\tau(y))^\sigma.
\]
Define the function $\mu: (0,1) \to (0,1)$ by
\[
\mu(\eta) = \min_{1 \leq i \leq q} \{ (\eta/|a_i|)^{1/i} \}.
\]
Then $\mu$ is a continuous monotone 
bijection, so it has an inverse function $\lambda$.
Since $\mu(\eta)> \eta$, we have $\lambda(\eta) < \eta$ for all $\eta$.
Moreover, for $\mu \in (0,1)$ sufficiently close to 1, the term
$i=q$ dominates and we have $\mu(\eta) = \eta^{1/q}$.

\begin{lemma} \label{lem:radius}
Suppose $\sigma$ is zero-centered. Then for any $x \in D$,
if $y$ is an element of $D$ of minimum norm for the property that
$\tau(y) = x$, then $\lambda(|y|) = |x|$.
\end{lemma}
\begin{proof}
Exercise with Newton polygons.
\end{proof}
\begin{cor}
If $a \in \Gcon$ and $a^\sigma$ converges for $\delta < |t| < 1$
for some $\delta > 0$, then $a$ converges for $\lambda(\delta) < |t| < 1$.
\end{cor}

\begin{lemma} \label{lem:64b}
If $a \in \Gcon$ and $f = a^\sigma/a$ is a rational function of $t$,
then $a$ factors as a rational function of $t$ times a unit of $\Omega$.
\end{lemma}
\begin{proof}
We may suppose $\sigma$ is zero-centered.
We will also allow tensoring over $W$ by an unspecified finite extension
of $W$ (which will be necessary in order to factor some polynomials);
we will use the same symbols $\Omega$ and $\Gcon$ for the resulting rings,
and extend $\sigma$ to these rings via the chosen extension of 
$\sigma$ to $W^{\alge}$.
Also, write $\ord(m,x)$ for the order of vanishing of an analytic function
$m$ at a point $x$ where it is defined.

Suppose the contrary, and choose a counterexample with $f$ of minimal
total degree. Then all zeroes and poles of $f$ must lie in $D$.
Write $f = g/h$, where $g$ and $h$ are relatively prime polynomials, so that
$a^\sigma h = ag$. 
Introduce an equivalence relation on $D$ defined by
$r \sim r'$ if $\tau(r) = \tau(r')$. Let $S_1, \dots, S_n$
be the equivalence classes containing roots of $h$.
For $i=1, \dots, n$,
let $r_{i1}, \dots, r_{i\ell}$ be the elements of $S_i$ (we should really
write $\ell_i$ for $\ell$, because $\ell$ depends on $i$, but we will suppress
that subscript for notational simplicity). For $j=1,\dots,\ell$,
let $c_{ij}$ and $m_{ij}$ be the multiplicities of $r_{ij}$ as roots of
the equations (in $x$) $\tau(x) = \tau(r_{i1})$ and $h(x) = 0$, respectively;
then $c_{ij} > 0$, $m_{ij} \geq 0$, $\sum_{j=1}^l c_{ij} = q$ (by Weierstrass
preparation), and
$\max_j \{ m_{ij}\} > 0$ (otherwise the class $S_i$ would not have been
labeled as such).
Put 
\[
m_i = \max_j \left\{ \left\lceil \frac{m_{ij}}{c_{ij}} \right\rceil \right\}.
\]

Let $e$ be the polynomial whose roots are $\tau(r_{i1})$ with
multiplicity $m_i$ for $i=1, \dots, n$. Then
$e^\sigma$ factors by Weierstrass preparation as a unit of 
$\Omega$ times a polynomial whose
roots are the $r_{ij}$, each with multiplicity $m_i$.
In particular,
$h$ divides the polynomial factor of $e^\sigma$.

Choose $\delta>0$ such that $a$ and $a^{-1}$ converge for $\delta < |t| < 1$.
We have 
\[
(ae)^\sigma = ag(e^\sigma/h),
\]
and the right side converges
for $\delta < |t| < 1$. Hence $ae$ converges for $\lambda(\delta)
< |t| < 1$.
Suppose $e$ has a zero $s$ with $\lambda(\delta) < |s| < 1$.
Then $s = \tau(r_{ij})$ for
some $i,j$, and we would have $\delta < |r_{ij}| < 1$ for all
such $j$ by Lemma~\ref{lem:radius}.

Choose $i$ and $j$ such that $s = \tau(r_{ij})$ and
$m_i = \lceil m_{ij}/c_{ij} \rceil$, so that in particular $m_{ij} > 0$.
Then $r_{ij}$ cannot be a root of $g$, since $g$ and $h$ are coprime 
polynomials and $r_{ij}$ is a root of $h$. Since $a$ and $a^{-1}$
both converge in the region $\delta < |t| < 1$ containing $r_{ij}$, we have
\begin{align*}
\ord((ae)^\sigma, r_{ij}) &= \ord(ag(e^\sigma/h),r_{ij}) \\
&= \ord(e^\sigma/h, r_{ij}) \\
&= m_i c_{ij} - m_{ij} \\
&\in \{0, \dots, c_{ij}-1\}.
\end{align*}
But $\ord((ae)^\sigma, r_{ij})$ is divisible by $c_{ij}$, so it must
be zero. Hence $(ae)^\sigma$ does not vanish at $r_{ij}$, so
$ae$ does not vanish at $\tau(r_{ij})$, and so $(ae)^\sigma$ does not
vanish at $r_{ij'}$ for $j' = 1,\dots,\ell$.
(Beware that we cannot conclude that $a$ does not vanish at $\tau(r_{ij})$
because it is not known to converge there!)
In particular, $e^\sigma/h$ does not vanish at $r_{ij'}$ for $j'=1,\dots,
\ell$; hence $m_i c_{ij'} = m_{ij'}$ for $j'=1,\dots,\ell$.

Now $((t-s)^{m_i})^\sigma$ factors by Weierstrass preparation as a
unit of $\Omega$ times a polynomial with roots $r_{i1}, \dots, r_{i\ell}$
of respective multiplicities $m_i c_{i1}, \cdots, m_i c_{i\ell}$.
The latter polynomial is equal to $h$ times a scalar,
because we just showed that $m_i c_{ij} = m_{ij}$ for $j=1, \dots, \ell$.
If we now write
\[
\frac{g((t-s)^{m_i})^\sigma}{h(t-s)^{m_i}} = 
\frac{(a(t-s)^{m_i})^\sigma}{a(t-s)^{m_i}},
\]
the left side is a unit of $\Omega$ times a rational function with
total degree strictly smaller than that of $f$. This yields a smaller 
counterexample to the original assertion than the one we chose,
contradicting that choice. Thus it cannot be that $e$ has a zero
$s$ with $\lambda(\delta) < |s| < 1$.

Since $ae$ converges for $\lambda(\delta) < |t| < 1$ and $e$ is
a polynomial with no zeroes in $\lambda(\delta) < |t| < \infty$,
$a$ must also converge for $\lambda(\delta) < |t| < 1$. By a
similar argument, $a^{-1}$ also converges for $\lambda(\delta) < |t| < 1$.
But we can now repeat the whole argument with $\delta$ replaced by
$\lambda(\delta)$, and so on. Since $\lambda$ is a continuous monotone
bijection
of $(0,1)$ onto itself and $\lambda(\eta) < \eta$ for all $\eta$,
the sequence $\delta, \lambda(\delta), \lambda(\lambda(\delta)), \cdots$
must converge to zero. In other words, $a$ and $a^{-1}$ converge
for $0 < |t| < 1$. By Lemma~\ref{lem:nonvanish},
$a = cut^n$ for some $c \in W$, some integer $n$ and some principal
unit $u$ of $\Omega$.
But then $a$ and $f$ did not form a counterexample to the desired
assertion to begin with, contradiction. Thus the claim of the lemma
holds, as desired.
\end{proof}

For our next lemma, we must move to a ring that relates to $\calR^+$
the same way that $\Gcon$ relates to $\Omega$. Namely,
let $\calR$ denote the ring of formal Laurent
series $\sum_{n\in \ZZ} c_n t^n$ with $c_n \in 
W [\frac 1p]$, such that $v(c_n) + rn \to \infty$ as $n \to -\infty$ for some
$r>0$ and $v(c_n) + sn \to \infty$ as $n \to +\infty$ for every
$s>0$. 
Each element
of $\calR$ can be viewed as a rigid analytic function on some unspecified
open annulus of outer radius 1.
The ring $\calR$ occurs commonly in the theory of $p$-adic
differential equations, where it is known as the Robba ring;
it contains both $\Gamma_c$ and $\calR^+$.

\begin{remark}
The conclusion of Lemma~\ref{lem:64b}
fails if we allow $a \in \calR$; for instance,
if $f = (1+t/p)^{-1}$, we can take $a$ to be the convergent product
$\prod_{i=1}^\infty (1+t/p)^{\sigma^i}$. 
\end{remark}

\begin{lemma} \label{lem:64a}
If $x,y \in \calR$ nonzero
satisfy $x^\sigma = gx$ and $y^\sigma = gp^{\ell} y$
for some $g \in \Gcon[\fp]$ 
and some nonnegative integer $\ell$,
then $\ell=0$ and $y/x \in \QQ_p^*$.
\end{lemma}
\begin{proof}
We may suppose $\sigma$ is zero-centered.
As in \cite{bib:lazard} or \cite{bib:mecrew},
we can factor $x = abb^\sigma b^{\sigma^2}\cdots$ for some
$a \in \Gcon[\frac 1p]$ 
and some $b \in K[t]$ with constant coefficient 1 having
all roots in the closed unit disc, as follows.
Pick $\delta < 1$ large enough so that
$\mu(\eta) = \eta^{1/q}$ for
$\eta \in (\delta,1)$,
and take the roots of $b$ to be the zeroes $r$ of $x$ (with
multiplicity) for which
$\eta < |r| \leq \eta^{1/q}$.
Then the roots of $b^{\sigma^i}$ in $D$ lie in the annulus
$\eta^{1/q^i} < |t| \leq \eta^{1/q^{i+1}}$ by
Lemma~\ref{lem:radius}; hence
$b^{\sigma^i}$
and $b^{\sigma^j}$ are relatively prime in $\Omega[\fp]$ whenever 
$i \neq j$.
Moreover, 
the infinite product $bb^\sigma b^{\sigma^2}\cdots $ converges in 
$\calR^+$ (see Remark~\ref{rem:product}), $x$ is divisible by
the limit (since $x$ is divisible by each factor and the factors
are pairwise coprime), 
and the quotient $a$ is a unit in $\calR$, hence is a nonzero
element of $\Gcon[\frac 1p]$.
Factor $y= cdd^\sigma d^{\sigma^2} \cdots$ analogously.

Now $g = a^\sigma/(ab)$ and $g p^{\ell} = c^\sigma/(cd)$, so we have an
equality
\[
p^{\ell} \frac{(a/c)^\sigma}{a/c} = b/d
\]
within $\Gcon[\frac 1p]$. 
By Lemma~\ref{lem:64b}, $a/c$ is a rational function of $t$ times
a unit of $\Omega$. Write $a/c = hu(r/s)$,
where $h \in W[\fp]$, $u$ is a principal unit in $\Omega$,
and
$r$ and $s$ are polynomials in $t$ with constant coefficient 1.
We then have the equality
\begin{equation} \label{eq:lem64}
bs^\sigma r = p^{\ell} (dr^\sigma s) \frac{(hu)^\sigma}{hu}
\end{equation}
in $\Omega[\fp]$; since $bs^\sigma r$ and $dr^\sigma s$ have
constant coefficient 1, so does $p^{\ell} (hu)^\sigma/(hu)$.
In other words, $p^{\ell} h^\sigma = h$, which implies
that $\ell = 0$ and $h \in \QQ_p^*$.

From \eqref{eq:lem64},
we have
\[
\prod_{i=0}^{n-1} (bs^\sigma r u)^{\sigma^i}
= \prod_{i=0}^{n-1} (dr^\sigma s u^\sigma)^{\sigma^i},
\]
which upon cancellation yields
\[
u s^{\sigma^n} r \prod_{i=0}^{n-1} b^{\sigma^i} = 
u^{\sigma^n} r^{\sigma^n} s \prod_{i=0}^{n-1} d^{\sigma^i};
\]
taking limits (thanks to Remark~\ref{rem:product} again) yields
\[
u r (x/a) = s (y/c).
\]
That is, $y/x = (cur)/(as) = 1/h \in \QQ_p^*$, as desired.
\end{proof}

\begin{remark} \label{rem:error}
The errors in \cite{bib:methesis} and \cite{bib:mz} referred to earlier
lie in their analogues of Lemma~\ref{lem:64a}.
In \cite[Lemma~3.2.4]{bib:methesis}, it is argued that
for $x,y \in \calR$, $(x/y)^\sigma = c(x/y)$ implies $x/y \in W$
by noting that this would be true if $x/y$ were a formal Laurent series.
However, there is no natural way to perform this division in a fashion
compatible with the action of $\sigma$.
In \cite[Lemma~32]{bib:mz}, there is an attempt to rectify this 
by clarifying that the division $x/y$ takes place in the ring $R_r$ of 
series $\sum c_n t^n$ with $v(c_n) + rn \to \infty$ as $n \to \pm \infty$
for some particular value of $r$. The subtlety here is that the
key equation $x^\sigma/y^\sigma = (x/y)^\sigma$ only holds when
the division on the left is performed in $R_r$ and the division
on the right is performed in $R_{r/p}$. Thus knowing that
$y^\sigma/y = x^\sigma/x$ does not imply that $(x/y)^\sigma = (x/y)$
in any meaningful sense: the symbols $x/y$ on the left and right
side are computed in different ways, so represent different
formal Laurent series.
\end{remark}

With Lemma~\ref{lem:64a} in hand, we can now establish Theorem~\ref{thm:64},
which will complete the proof of Theorem~\ref{thm:main3}.
\begin{proof}[Proof of Theorem~\ref{thm:64}]
We begin with some reductions as 
in \cite[Proposition~6.4]{bib:dej4}.
We first reduce to the case of $M$ free.
By Lemma~\ref{lem:free}, there is an isogeny $\psi: M \to M'$ with $M'$ free
as an $\Omega$-module. Let $\psi_c: M_c \to M'_c = M' \otimes \Gcon$
be the induced map, and let $N'_c$ be the saturation of $\psi_c(N_c)$ in 
$M'_c$. If Theorem~\ref{thm:64} is already known for free modules, then
$N'_c = N' \otimes \Gcon$ for some saturated sub-$F$-module $N'$ of $M'$.
Let $N$ be the saturation of $\psi^{-1}(N')$; then $N \otimes \Gcon = N_c$,
as desired.

We may thus assume that $M$ is free. We next reduce to the case where
$N_c$ has rank 1. Suppose instead that $\rank N_c = d > 1$.
If the rank 1 case is already known, we may apply it to
$\wedge^d N_c \subseteq \wedge^d M_c = (\wedge^d M) \otimes \Gcon$
to obtain a saturated sub-$F$-module $N_1$ of $\wedge^d M$ such that
$N_1 \otimes \Gcon = \wedge^d N_c$.
Let $K$ be the $p$-adic valuation subring of $\Frac \Omega$.
The set of $\bv \in M_c$ such that $\bv \wedge \bw = 0$ for
any $\bw \in N_1$ is exactly $N_c$; this set is defined by linear conditions,
so the set of $\bv \in M \otimes K$ such that $\bv \wedge \bw = 0$
for any $\bw \in N_1$ spans $N_c$ over $\Gcon$. In particular, if we
let $N$ denote the set of $\bv \in M$ such that $\bv \wedge \bw = 0$
for any $\bw \in N_1$, then $N$ is a saturated sub-$F$-module of $M$
which spans $N_c$ over $\Gcon$. Hence $N \otimes \Gcon = N_c$, as desired.

We may thus reduce to the case where $M$ is free and $N_c$ has rank 1.
By Proposition~\ref{prop:dwork},
we can find a basis $\bv_1, \dots, \bv_n$ of $M \otimes \calR^+$
such that $F^m\bv_i = p^{\ell_i} \bv_i$ for some integers $m>0$ and $\ell_i
\geq 0$.
Choose a generator $\bw$ of $N_c$; in $M \otimes \calR$, we then
have the equality
$\bw = \sum_i h_i \bv_i$ for some $h_i \in \calR$. 
Assume without loss
of generality that $h_1 \neq 0$.
Then $F^m \bw = g\bw$ for some nonzero $g 
\in \Gcon$, which implies
\begin{align*}
g\bw &= F^m \bw \\
&= \sum_i F^m (h_i \bv_i) \\
&= \sum_i h_i^{\sigma^m} p^{\ell_i} \bv_i.
\end{align*}
That is, $h_i^{\sigma^m} = g p^{-\ell_i} h_i$ for each $i$.
By Lemma~\ref{lem:64a} applied to $x = h_i$ and $y = h_1$,
we see that if $h_i \neq 0$, then $\ell_1 = \ell_i$ and
$h_i/h_1 \in \QQ_p^*$. In particular, each
$h_i$ is divisible by $h_1$; on the other hand,
the $h_i$ generate the unit ideal
since $\bw$ is part of a basis of $M \otimes \calR$.
(That is because $\bw$ is already part of a basis of $M_c$,
which is true because $\bw$ generates a saturated submodule of $M_c$
and $\Gcon$ is a principal ideal domain.)
Hence $h_1$ is a unit in $\calR$, that is, $h_1 \in \Gcon[\frac 1p]$.

Now set $\bv = (p^r/h_1) \bw = 
\sum_i p^r (h_i/h_1) \bv_i$, where $r$ is the smallest
integer for which the quantity on the right lies in $M_c$
(rather than $M \otimes \Gcon[\fp]$).
Then $F^m \bv$ is a multiple
of $\bv$, but by comparing the coefficients of $\bv_1$ in $\bv$ and
$F^m\bv$, we see that in fact $F^m \bv = p^{\ell_1} \bv$.
Thus $\bv$ is in the $\ZZ_p$-span of the $\bv_i$, and so belongs to
$M \otimes_\Omega (\calR^+ \cap \Gcon ) 
= M \otimes_\Omega \Omega
= M$. Let $N$ be the $\Omega$-span of $\bv$ in $M$;
then $N$ is saturated and $N_c = N \otimes \Gcon$, as desired.
\end{proof}

We conclude this section with two examples illustrating subtleties
in the above arguments. First, observe that
for $\sigma$ standard and $n$ an integer, 
the equation $x^\sigma = p^n x$ has a solution $x \in \calR$
if and only if $n = 0$.
However, consider the Frobenius lift $\sigma$ defined by
\[
t^\sigma = (t+1)^p - 1,
\]
which arises naturally in the theory of $(\phi, \Gamma)$-modules
(see for instance \cite{bib:berger}). Then 
\[
x = \log (1+t) = \sum_{n=1}^\infty \frac{(-1)^{n-1} t^n}{n}
\]
clearly belongs to $\calR^+$ and $x^\sigma = px$. In particular, the
non-uniqueness visible in the proof of Proposition~\ref{prop:dwork}
(in the choices of the $V_h$ for $h < h_0$)
is more than just an artifact of the proof technique; the matrix $U$
therein really may fail to be unique.

Second, note that even for $\sigma$ standard, so that
$x^\sigma = p^\ell x$ has no solutions with $x \in \calR$ when
$\ell \neq 0$, it is still possible to have
$x^\sigma = gx$ for $x \in \calR$ and $g \in \Gcon$ of nonzero valuation.
For instance, if $x = bb^\sigma b^{\sigma^2} \cdots$ with $b = (1+t/p)$,
then $x^\sigma = gx$ for
\[
g = (1+t/p)^{-1} = p/t (1+p/t)^{-1} = \sum_{i=0}^\infty (-1)^i (p/t)^{i+1}.
\]

\section{A Witt vector counterexample}
\label{sec:witt}

Instead of considering $F$-modules, which depend on the choice of
a Frobenius lift, one might prefer to consider
modules over a ring that lifts $k\llbracket t \rrbracket$ and comes
equipped with a \emph{canonical} Frobenius. In crystalline cohomology, 
this point of view leads to the de Rham-Witt complex introduced by
Deligne-Illusie \cite{bib:illusie}. One has an analogous construction
in our context: given an $F$-crystal over $\Omega$, one obtains
by base extension a finite free $W(k \llbracket t \rrbracket)$-module
$M$ equipped with an isogeny
$F: \sigma^* M \to M$, where $\sigma$ here denotes the Witt vector
Frobenius, and the isogeny property means there exists an
additive map $V: M \to \sigma^* M$ such that $F \circ V = V \circ F$
is multiplication by some power of $p$. 

From the $F$-crystal structure, one can also obtain a ``de Rham-Witt
connection'' that maps $M$ to the tensor product of $M$ with
a suitable quotient of the module of K\"ahler differentials of
$W(k \llbracket t \rrbracket)$, on which the Frobenius and Verschiebung
operators of the Witt ring act.
One can use this extra data to give an equivalent formulation of de 
Jong's theorem that does not require an auxiliary choice of Frobenius
lift; we will not work this out here
(but compare the construction of Gauss-Manin connections in the de Rham-Witt
context given by Langer and Zink in \cite{bib:langer zink}).

However, as in the introduction, one can find situations in this context
in which a Frobenius structure naturally arises by itself, e.g., in
Zink's theory of displays of $p$-divisible groups \cite{bib:zink}.
One can then ask whether one has an analogue of
Theorem~\ref{thm:main2} to the effect that given a free $F$-module
$M$ over $W(k \llbracket t \rrbracket)$ and
an element $\bv \in M \otimes W(k((t)))$ 
such that $F\bv = p^{\ell} \bv$ for some
integer $\ell$, it follows that $\bv \in M$.
This would in fact generalize Theorem~\ref{thm:main2} if it were true;
however, we can exhibit a counterexample as follows.

\begin{prop} \label{prop:counter}
There exists a free $W(k \llbracket t \rrbracket)$-module $M$ of rank $2$,
an isogeny $F: \sigma^* M \to M$ and an element $\bv$ of 
$(M \otimes W(k ((t)))) \setminus M$ such that $F\bv = \bv$.
\end{prop}
\begin{proof}
Put $R = W(k \llbracket t \rrbracket)$ for brevity, let brackets denote
Teichm\"uller lifts, and let $V$ denote the Verschiebung map on $R$.
We first define sequences $\{x_n\}_{n=0}^\infty$ and
$\{b_n\}_{n=0}^\infty$ over $R$
satisfying
\begin{align*}
b_n &\equiv [t] \pmod{R^V}, \\
x_n &\equiv [t]^{p^2-p} \pmod{pR}, \\
-b_n^{\sigma^2} + x_n b_n^\sigma - p^2 b_n &\equiv 0 \pmod{p^2 R^{V^n}},
\end{align*}
as follows.
Begin by setting $x_0 = [t]^{p^2-p}$ and $b_0 = [t]$.
Given $x_n$ and $b_n$, define $\Delta_n \in R$ by the equation
$p^2 \Delta_n^{V^n} = 
-b_n^{\sigma^2} + x_n b_n^\sigma - p^2 b_n$,
and set
\[
b_{n+1} = b_n + \Delta_n^{V^{n+2}}, \qquad
x_{n+1} = x_n - [t]^{p^2 - 2p} p \Delta_n^{V^{n+1}}.
\]
We then have
\begin{align*}
-b_{n+1}^{\sigma^2} + x_{n+1} b_{n+1}^\sigma - p^2 b_{n+1}
&\equiv -b_n^{\sigma^2} - \Delta_n^{V^{n+2} \sigma^2} + 
x_{n+1} b_{n+1}^\sigma - p^2 b_n - p^2 \Delta_n^{V^{n+2}} \\
&\equiv -b_n^{\sigma^2} + x_n b_n^\sigma - p^2 b_n
- \Delta_n^{V^{n+2} \sigma^2} + 
x_{n+1} b_{n+1}^\sigma - x_n b_n^\sigma \\
&\equiv p^2 \Delta_n^{V^n} - p^2 \Delta_n^{V^n}
+ x_{n+1} b_{n+1}^\sigma - x_n b_n^\sigma \\
&\equiv (x_{n+1} - x_n) b_{n+1}^\sigma
+ x_n(b_{n+1} - b_n)^\sigma \\
&\equiv -[t]^{p^2-2p} 
p\Delta_n^{V^{n+1}} b_{n+1}^\sigma + p x_n \Delta_n^{V^{n+1}} \\
&\equiv p \Delta_n^{V^{n+1}} (-[t]^{p^2-2p} b_{n+1}^\sigma + x_n) \\
&\equiv 0 \pmod{p^2 R^{V^{n+1}}},
\end{align*}
the last congruence holding because $[t]^{p^2-2p} b_{n+1}^\sigma 
\equiv [t]^{p^2-p} \pmod{pR}$ and $x_n \equiv [t]^{p^2-p}
\pmod{pR}$.

Note that in the above construction, $x_{n+1} \equiv x_n \pmod{p R^{V^{n+1}}}$
and $b_{n+1} \equiv b_n \pmod{R^{V^{n+2}}}$. Thus
the sequences $\{x_n\}$ and $\{b_n\}$ converge in $R$ to
limits $x$ and $b$, respectively, satisfying
\[
-b^{\sigma^2} + x b^\sigma = p^2 b.
\]
We now define the module $M$ to be the module of rank 2 column
vectors over $R$ with the $\sigma$-linear map $F$ given by
\[
F \begin{pmatrix} y \\ z \end{pmatrix} =
\begin{pmatrix} 0 & p^2 \\ -1 & x \end{pmatrix}
\begin{pmatrix} y^\sigma \\ z^\sigma \end{pmatrix}.
\]
By construction, the column vector
\[
\bw = \begin{pmatrix} b^\sigma \\ b \end{pmatrix} 
\]
then satisfies $F \bw = p^2 \bw$. 

We next construct a sequence $\{d_n\}_{n=1}^\infty$ of elements of 
$S = W(k((t))^{\sep})$ such that $d_{n+1} - d_n \in S^{V^n}$ and
\[
-p^2 d_n^{\sigma^2} + x d_n^\sigma - d_n \equiv 0 \pmod{S^{V^n}}.
\]
To begin, set $d_1 = [t]^{-p}$.
Given $d_n$, define $\Delta_n$ by the equation
$\Delta_n^{V^n} = -p^2 d_n^{\sigma^2} + x d_n^\sigma - d_n$.
Since the polynomial $z^p - z - c$ in $z$ is separable over $k((t))^{\sep}$
for any $c \in k((t))^{\sep}$,
we can find $z_n \in S$ such that
$z_n^p - z_n + [t]^{p^{n+1}} \Delta_n \equiv 0 \pmod{S^V}$.
Now put $d_{n+1} = d_n + [t]^{-p} z_n^{V^n}$; then
\begin{align*}
-p^2 d_{n+1}^{\sigma^2} + x d_{n+1}^\sigma - d_{n+1}
&\equiv
\Delta_n^{V^n} - p^2 ([t]^{-p} z_n^{V^n})^{\sigma^2}
+ x ([t]^{-p} z_n^{V^n})^\sigma - [t]^{-p} z_n^{V^n}
\\ &\equiv
\Delta_n^{V^n} + x [t]^{-p^2} z_n^{V^n \sigma} - [t]^{-p} z_n^{V^n} \\
&\equiv \Delta_n^{V^n} + p [t]^{-p} z_n^{V^{n-1}} - [t]^{-p} z_n^{V^n} \\
&\equiv (\Delta_n^V + p [t]^{-p^n} z_n - [t]^{-p^n} z_n^V)^{V^{n-1}} \\
&\equiv (\Delta_n + [t]^{-p^{n+1}} z_n^p - [t]^{-p^{n+1}} z_n)^{V^{n}} \\
&\equiv 0 \pmod{S^{V^{n+1}}}.
\end{align*}
The sequence $\{d_n\}$ converges to a limit $d \in S$ such that
\[
-p^2 d^{\sigma^2} + x d^\sigma = d.
\]
Moreover, any $D \in S$ such that
$-p^2 D^{\sigma^2} + x D^\sigma = D$ must be of the form
$dr$ for some $r \in \ZZ_p$. Namely, it suffices to check that
if $D \in S^{V^n}$ satisfies the equation, then $D - rd \in
S^{V^{n+1}}$ for some $r \in \ZZ_p$; but the set of all such
$D$ generates a subgroup of $S^{V^n}/S^{V^{n+1}}$ of order $p$,
as do the multiples of $p^n d$.

Now put $e = b^\sigma d - p^2 bd^\sigma \in S$; then
\begin{align*}
e^\sigma  &=
b^{\sigma^2} d^\sigma - p^2 b^\sigma d^{\sigma^2} \\
&=  d^\sigma( x b^\sigma - p^2 b) - b^\sigma (xd^\sigma - d) \\
&= b^\sigma d - p^2 bd^\sigma = e,
\end{align*}
and so $e \in \ZZ_p$.
Given $\tau \in \Gal(k((t))^{\sep}/k((t)))$, we have a natural
action of $\tau$ on $S$ commuting with $\sigma$, and so
$-p^2 d^{\tau \sigma^2} + x d^{\tau \sigma} = d^\tau$.
As noted above, this implies that $d^\tau = \chi(\tau) d$ for some
$\chi(\tau) \in \ZZ_p$; in fact, $\chi(\tau) \in \ZZ_p^*$ since
$d$ and $d^\tau$ are both nonzero modulo $S^V$. 
However,
\begin{align*}
e &= e^\tau \\
&= (b^\sigma d - p^2 bd^\sigma)^\tau \\
&= b^\sigma d^\tau - p^2 b d^{\sigma \tau} \\
&= \chi(\tau) (b^\sigma d - p^2 bd^\sigma) \\
&= \chi(\tau) e,
\end{align*}
forcing $\chi(\tau) = 1$ for all $\tau$. That is,
$d$ actually belongs to the subring of $S$ fixed by
$\Gal(k((t))^{\sep}/k((t)))$, which is precisely
$W(k((t)))$.

To conclude, note that the column vector
\[
\bv = \begin{pmatrix} p^2 d^\sigma \\ d \end{pmatrix}
\]
of $M \otimes W(k((t)))$ satisfies $F\bv = \bv$, but it is evident
that $\bv \notin M$ since $d \equiv [t]^{-p} \pmod{W(k((t)))^V}$.
This yields the desired counterexample.
\end{proof}

We note in concluding that the counterexample constructed above
is not quite as magical as it might look. In the language of
\cite{bib:mecrew}, $\bv$ and $\bw$
generate the first steps in the descending and ascending slope filtrations;
all that was necessary was to adjust $M$ so that $\bv$ would come out
defined over $W(k \llbracket t \rrbracket)$, not just over
$W(k \llbracket t \rrbracket^{\alge})$. Moreover, the ``generic'' slopes
had to be taken more than 1 apart, otherwise $M$ would automatically
admit a ``de Rham-Witt connection'' and de Jong's theorem would apply; 
indeed, this is what happens in
the case of displays of $p$-divisible groups, in which the slopes are
constrained to lie in $[0,1]$.


\begin{thebibliography}{AdV}
\bibitem[AdV]{bib:andre divizio}
Y. Andr\'e and L. di Vizio, $q$-difference equations and $p$-adic local
monodromy, to appear in \emph{Ast\'erisque}; preprint available at
\texttt{picard.ups-tlse.fr/\~{}divizio}.

\bibitem[B]{bib:berger}
L. Berger, Repr\'esentations $p$-adiques et \'equations diff\'erentielles,
\textit{Invent. Math.} \textbf{148} (2002), 219--284.

\bibitem[dJ]{bib:dej4}
A.J. de~Jong, Homomorphisms of Barsotti-Tate groups and crystals in
positive characteristic, \textit{Invent. Math.} \textbf{134} (1998),
301--333; erratum, \textit{ibid.} \textbf{138} (1999), 225.

\bibitem[I]{bib:illusie}
L. Illusie,
Complexe de de Rham-Witt et cohomologie cristalline,
\textit{Ann. Sci. \'Ec. Norm. Sup.} \textbf{12} (1979), 501--661.

\bibitem[Ka]{bib:katz}
N.M. Katz, Slope filtrations of $F$-crystals,
\textit{Ast\'erisque} \textbf{63} (1979), 113--163.

\bibitem[Ke1]{bib:methesis}
K.S. Kedlaya, Descent theorems for overconvergent $F$-crystals,
Ph.D. thesis, Massachusetts Institute of Technology, 2000;
available at \texttt{math.mit.edu/\~{}kedlaya}.

\bibitem[Ke2]{bib:mecrew}
K.S. Kedlaya, A $p$-adic local monodromy theorem, to appear in
\textit{Ann. of Math.}; \texttt{arXiv:
math.AG/0110124}.

\bibitem[Ke3]{bib:mefull}
K.S. Kedlaya, Full faithfulness for overconvergent $F$-isocrystals,
in A. Adolphson et al. (eds.),
\textit{Geometric Aspects of Dwork Theory (Volume II)},
de Gruyter (Berlin), 2004, 819--835.

\bibitem[LZ]{bib:langer zink}
A. Langer and Th.\ Zink,
De Rham-Witt cohomology for a proper and smooth morphism,
\textit{J. Inst. Math. Jussieu} \textbf{3} (2004), 231--314.

\bibitem[La]{bib:lazard}
M. Lazard, Les z\'eros des fonctions analytiques d'une variable sur
un corps valu\'e complet, \textit{Publ. Math. IHES} \textbf{14} (1962), 47--75.

\bibitem[Ma]{bib:manin}
Yu.I. Manin, The theory of commutative formal groups over fields of
positive characteristic, \textit{Russian Math. Surveys} \textbf{18}
(1963), 1--83.

\bibitem[MZ]{bib:mz}
W. Messing and Th.\ Zink, de Jong's theorem on homomorphisms of
$p$-divisible groups, preprint, 2001; available at
\texttt{www.mathematik.uni-bielefeld.de/\~{}zink}.

\bibitem[Z]{bib:zink}
Th.\ Zink, The display of a formal $p$-divisible group,
in Cohomologies $p$-adiques et applications arithm\'etiques, I,
\textit{Ast\'erisque} \textbf{278} (2002), 127--248.

\end{thebibliography}
\end{document}